\documentclass{amsart}
\usepackage{amsfonts,amssymb,amscd,amsmath,enumerate,verbatim,calc}

\newcommand{\IFF}{\text{if and only if}}
\newcommand{\wrt}{with respect to}

\newcommand{\n}{\mathfrak{n} }
\newcommand{\m}{\mathfrak{m} }
\newcommand{\M}{\mathfrak{M} }

 \newcommand{\rt}{\rightarrow}

\newcommand{\ov}{\overline}

\newcommand{\wt}{\widetilde }

\newcommand{\reg}{\operatorname{reg}}

\newcommand{\T}{\text}

\newcommand{\depth}{\operatorname{depth}}

\newcommand{\pd}{\operatorname{projdim}}
\newcommand{\G}{\operatorname{G-dim}}
\newcommand{\Hd}{\operatorname{H-dim}}
\newcommand{\CI}{\operatorname{CI-dim}}
\newcommand{\CIt}{\operatorname{CI_{*}-dim}}
\newcommand{\Gt}{\operatorname{G^*-dim}}
\newcommand{\Ct}{\operatorname{CM-dim}}
\newcommand{\cx}{\operatorname{cx}}

\theoremstyle{plain}
\newtheorem{theorem}{Theorem}[section]
\newtheorem{corollary}[theorem]{Corollary}
\newtheorem{lemma}[theorem]{Lemma}
\newtheorem{proposition}[theorem]{Proposition}

\newtheorem{thm}{Theorem}

\theoremstyle{definition}
\newtheorem{definition}[theorem]{Definition}
\newtheorem{s}[theorem]{}
\newtheorem{remark}[theorem]{Remark}

\theoremstyle{remark}

\numberwithin{equation}{section}
\begin{document}

\title{An analogue of a theorem due to Levin and Vasconcelos}
\author{Javad Asadollahi, Tony ~J.~Puthenpurakal}
\date{\today}

\address{School of Mathematics, Institute for studies in Theoretical Physics
and Mathematics (IPM), P.O.Box: 19395-5746,Tehran, Iran and
Shahre-Kord University, P.O.Box: 115, Shahre-Kord, Iran.}
\email{Asadollahi@ipm.ir}

\address{Department of Mathematics, IIT Bombay,
Powai , Mumbai 400076.}
\email{tputhen@math.iitb.ac.in}

\keywords{Homological dimensions, Gorenstein dimension, CI dimension, CM-dimension, $\m$-full ideal, superficial elements, Ratliff-Rush filtration }
\subjclass{ Primary 13H05, Secondary 13D10}
\begin{abstract} Let $(R,\m)$ be a Noetherian local ring. Consider
the notion of homological dimension of a module, denoted H-dim,
for H= Reg, CI, CI$_*$, G, G$^*$ or CM. We prove that, if for a finite
$R$-module $M$ of positive depth, $\Hd_R({\m}^iM)$ is finite for
some $i \geq \reg(M)$, then the ring $R$ has property H.
\end{abstract}
\maketitle
\section*{introduction}
One of the most influential results in commutative algebra is
 result of Auslander, Buchsbaum and Serre: A local ring $(R, \m)$ is regular if and only if
$\pd_RR/\m$ is finite. This result is considerably strengthened by
Levin and Vasconcelos. Their result can be read as follows,

\noindent \textbf{Theorem}  (Levin-Vasconcelos \cite{LV}).
\textit{Let $(R, \m)$ be a commutative Noetherian local ring. If
there exists a finite module $M$ such that $\pd_R {\m}^iM $ is
finite, for some $i\geq 1$, then $R$ is regular.}

Note that $\pd_R 0$ is defined to be $-\infty$. This in particular
implies that if $\dim R > 0$ then $R$ is regular if and only if
$\pd_R(R/{\m}^i)$ is finite for some $i \geq 1$.

On the other hand a variety of refinements of projective dimension
are defined, namely G-dimension \cite{AB}, CI-dimension
\cite{AGP}, CI$_*$-dimension \cite{Ger}, G$^*$-dimension
\cite{Vel} and CM-dimension \cite{Ger}. For the sake of uniformity
of notation sometimes we write Reg-dimension for projective
dimension. We say that $R$ has property H with H=Reg (respectively,
CI, CI$ _*$, G$ ^*$, G or CM) if it is regular (respectively,
complete intersection, complete intersection, Gorenstein,
Gorenstein or Cohen-Macaulay).  All the above homological dimensions
generalise   Auslander, Buchsbaum and Serre's result namely, the ring $R$ has
property H if and only if $\Hd_R R/\m$ is finite.

So a natural question arises, whether it is possible to
generalize Levin and Vasconcelos result to the above mentioned
generalizations of projective dimension. We prove the result when $\depth M > 0$
and $n \gg 0$. To state the result we need an invariant $\rho(M)$, see \ref{rr}.
In \cite{Pu3} it is proved that $\rho(M) \leq \reg(M)$.
Our main Theorem can be stated as follows.

\begin{thm}
\label{lev} Let $(R, \m)$ be a Noetherian local ring with residue
field $k$ and let $M$ be a finite $R$-module of positive depth.
Let H equal to  Reg (respectively, CI,CI$_*$, G$^*$,G, or CM).
 If
$\Hd_R(\m^nM)< \infty$ for some $n \geq \rho(M)$, then $R$ has
property H.
\end{thm}

As a corollary we get that if $\depth R > 0$ and $\Hd_RR/\m^n$ is
finite for some $n \geq  \reg(R)$, then $R$ has property H (see
Corollary \ref{Lcor}).

Recently Goto and Hayasaka \cite{GH} have proved that if
G-dimension of an integrally closed $\m$-primary ideal $I$ is
finite then $R$ is Gorenstein. In \cite{Pu4} Iyengar and the second author
used techniques used in the study of Hilbert functions to study some
homological properties of $M/\m^nM$ for $n \gg 0$. In this paper
we combine some of the techniques used in  these papers.

Here is an overview of the contents of the paper. In section 1 we introduce
notation and discuss a few preliminary facts that we need.
In Section two we define the notion of $\m$-full submodules, which
is a natural generalization of $\m$-full ideal, defined by D. Rees
(unpublished). See J. Watanabe's papers \cite{W1} and \cite{W2}
for some properties of $\m$-full ideals. In Section 3 we state the properties of H-dim
that we need and  prove a few that has not been explicitly proved before.
In Section 4 we prove the Theorem 1 and a few corollaries.

\section{Preliminaries}
All rings are commutative Noetherian and all modules are finite.
In this section we will give some auxiliary results that we need
for the proof of Theorem 1. Let $(R,\m)$ be a local ring with
residue field $k$ and let $M$ be a finite $R$-module. Let $G(R) = \bigoplus_{n \geq 0}
\m^n/\m^{n+1}$ be the associated graded ring of $R$ \wrt \ $\m$
and let $G(M) = \bigoplus_{n \geq 0} \m^nM/\m^{n+1}M$ be the
associated graded module of $M$ \wrt \ $\m$. Set $\M =
\bigoplus_{n \geq 1}\m^n/\m^{n+1}$ the maximal irrelevant ideal of
$G(R)$.

\s \textbf{Superficial elements:} An element $x \in \m$ is said to be $M$-\emph{superficial}
if there exists an integer $c > 0$ such that
$$ (\m^{n}M \colon_M x)\cap\m^cM =  \m^{n-1}M \ \text{ for all }\quad n > c
$$ When $\dim M >0$ it is easy to see that $x$ is superficial if and
only if $x\notin \m^2$ and $x^*$ does not belong to any relevant
associated prime of $G(M)$. Therefore, superficial elements always
exist if the residue field  $k$ is infinite. If $\depth M > 0$ and
$x$ is $M$-superficial then $(\m^{n}M \colon_M x)= \m^{n-1}M$ for
all $n \gg 0$ (see \cite[p.\ 7]{Sa} for $M =R$, the  general case
is similar).

\s \textbf{Regularity:} For $i \geq 0$ set $H^{i}_{\M}(G(M))$ to be the $i$'th
local cohomology of $G(M)$ \wrt \ $\M$. The modules
$H^{i}_{\M}(G(M))$ are graded and Artininan. Define
\[
\reg(M) = \max \{ i + j \mid H^{i}_{\M}(G(M))_j \neq 0 \}.
\]

\s If the residue field of $R$ is finite then we resort to the
standard trick to replace $R$ by $R' = R[X]_S$ where $S =
R[X]\setminus \m R[X]$. The maximal ideal of $R'$ is $\n = \m R'$.
The residue field of $R'$ is $l = k(X)$, the field of rational
functions over $k$. Set  $M'= M \otimes_R R'$. One can easily show
that $\n^i = \m^iR'$ for all $i \geq 1$, $\dim M = \dim M'$,
$\depth M = \depth M'$ and $\reg(M) = \reg(M')$.

\s \textbf{Ratliff-Rush Filtration }For every $n$ we consider the
chain of submodules of $M$
\[
\m^nM \subseteq \m^{n+1}M \colon_M \m \subseteq \m^{n+2}M \colon_M
\m^2 \subseteq \cdots \subseteq \m^{n+k}M \colon_M \m^k \subseteq
\cdots
\]
This chain stabilizes at a submodule which we denote by
\[
\wt{\m^n M} = \bigcup_{k \geq 1} \m^{n+k}M \colon_M \m^k.
\]
Clearly $\m \wt{\m^n M} \subseteq \wt{\m^{n+1} M}$ and
 $\wt{\m^{n+1} M} \subseteq \wt{\m^{n} M}$. So
$\mathcal{F} = \{ \wt{\m^n M} \} $ is an $\m$-filtration of $M$.
It is called the Ratliff-Rush filtration of $M$ \wrt \  $\m$.

\s   If $\depth M > 0$ then the following holds (see \cite{RR} for
the case $M = R$,  see \cite{Pu3} where it is proved in general)
\begin{enumerate}[\rm \quad 1.]
\item
$\wt{\m^nM} = \m^nM$ for all $n \gg 0$.
\item
$ \depth G(M) > 0$ \IFF \ $\wt{\m^nM} = \m^nM$ for all $n \geq 0$.

\item
If $x$ is  $M$-superficial  then

$(\wt{\m^{n+1}M}\colon_M x) = \wt{\m^nM}$ for all $n \geq 1$.

\end{enumerate}

\s
\label{rr}
In view of property 1. and 2. above it is convenient to define
the following invariant of $M$
\begin{equation*}
\rho(M) = \min \{ i \mid \wt{\m^nM} = \m^nM \ \ \text{for all} \ n
\geq i \}.
\end{equation*}
In \cite[Theorem 5]{Pu3} it is proved that $\rho(M) \leq \reg(M)$.

\section{ $\m$-full submodules}
In this section we generalize to modules the notion of $\m$-full
ideals. Throughout this section $R$ is local and $M$ is a finite
$R$-module.
\begin{definition}
A submodule $N$ of $M$  is called $\m$-full if there exists $x \in
\m$ such that $ \m N \colon_M x = N$.
\end{definition}
\begin{proposition}
\label{rrmfull} Let $(R,\m)$ be local with infinite residue field
and let $M$ be an $R$-module of positive depth. Then for all $n$
the submodule $\wt{\m^n M}$ is $\m$-full.
\end{proposition}
\begin{proof}
Since the residue field is infinite we can choose  $x \in \m$
which is  $M$-superficial \wrt \ $\m$. For all $n$ we have $\left(
\wt{\m^{n+1}M}\colon_M x \right) = \wt{\m^nM}$. So  we have
\[
  \wt{\m^n M} \subseteq \left(\m\wt{\m^{n}M}\colon_M x \right) \subseteq
\left( \wt{\m^{n+1}M}\colon_M x \right) =  \wt{\m^n M}
\]
So  $ \wt{\m^n M} $ is $\m$-full.
\end{proof}
In the next proposition we collect the basic properties of
$\m$-full submodules.
\begin{proposition}
\label{mfullp} Let $(R,\m, k)$ be local  , $M$ an $R$-module and
$N$ an $\m$-full submodule of $M$. Let $x \in \m$ be such that $\m
N \colon_M x = N$.  We have the following:
\begin{enumerate}[\rm (1)]
\item
$N \colon_M x = N \colon_M \m $.
\item
If $p_1,\ldots, p_l \in M$ is such that
$\{\ov{p_1},\ldots,\ov{p_l} \}$ is a basis of the $k$-vector space
$\left(N \colon_M \m \right)/N$  then $xp_1,\ldots, xp_l$ form
part of a minimal basis of $N$.
\item
Let $\{xp_1,\ldots, xp_l \}\bigcup \{z_1,\ldots,z_m \}$ generate
$N$ \emph{minimally}.

Define $\phi \colon N \rt \left(N \colon_M
\m \right)/N $ as follows: For $t = \sum_{i =1}^{l} a_i xp_i +
\sum_{j =1}^{m} b_j z_j$, set $\phi(t) = \sum_{i = 1}^{l}\ov{a_i}
\ov{p_i}$. The map $\phi$ is well defined and $R$-linear.
\item
The $R$-linear map $ \psi \colon \left(N \colon_M \m \right)/N \rt
N/xN $ defined by $ \psi(s + N) = xs + xN$ is a split injection.
\item
$k$ is a direct summand of $N/xN$
\end{enumerate}
\end{proposition}
\begin{proof}
(1) Clearly $N \colon_M x \supseteq N \colon_M \m $. Note that
\[
N \colon_M \m = \left(\m N \colon_M x \right) \colon_M \m =
\left(\m N \colon_M \m \right) \colon_M x \supseteq N\colon_M x.
\]
So  we get that $N \colon_M x = N \colon_M \m $.

(2)  If $\sum_{i =1}^{l} a_ixp_i \in \m N$ then $\sum_{i =1}^{l}
a_ip_i \in \left(\m N \colon_M x \right) = N$. Since
$\{\ov{p_1},\ldots,\ov{p_l} \}$ is a $k$-basis of $\left(N
\colon_M \m \right)/N$ we get that $a_i \in \m$ for all $i$.

(3) First note that if $t = \sum_{i =1}^{l} a_i xp_i + \sum_{j
=1}^{m} b_j z_j $, set $u = \sum_{i =1}^{l} a_i p_i$. Since
$xu = \sum_{i =1}^{l} a_i (xp_i)   \in N$ we have that $u \in N\colon_M x =
N\colon_M \m$.  Also note that if
 $t = \sum_{i =1}^{l} a_i xp_i + \sum_{j
=1}^{m} b_j z_j  = \sum_{i =1}^{l} a_{i}' xp_i + \sum_{j =1}^{m}
b_j' z_j$ then $a_i - a_i' \in \m $ for all $i$ and so $\sum_{i =
1}^{l}\left(\ov{a_i}  - \ov{a_i'} \right)\ov{p_i}  \in N $, since $p_i \in (N \colon_M \m )$ for all $i = 1,\ldots, l$. So it follows that
$\phi$ is well defined. Clearly $\phi $ is $R$-linear.

(4) Set $W = \left(N \colon_M \m \right)/N $.  Clearly $\psi $ is
$R$-linear. Furthermore by (2) we have that $\psi$ is injective.
Note that the map $\phi \colon N \rt W $ defined in (3) maps $xN$
to $0$ and so defines a map $\ov{\phi} \colon N/xN \rt W$.
Furthermore $\ov{\phi} \psi = id_W$. So $\psi$ is split.

(5) Since $\left(N \colon_M \m \right)/N $ is a $k$-vector space,
the assertion  follows from (4).
\end{proof}

\section{some properties of H-dim}
Let $(R,\m)$ be a local ring and let $M$ be a $R$-module.
For all unexpalined terminology see the
 survey paper \cite{Av}.

\begin{remark}
\label{HremL}
 It suffices to prove Theorem \ref{lev} when H is equal
to Reg (respectively, CI$_*$,  G, or CM) since homological dimensions satisfy the inequalities
\[
\Ct_R M \leq \G_R M \leq \CIt_R M \leq \CI_R M \leq \pd_R M
\]
\[
 \G_R M \leq \Gt_R M \leq \CI_R M
\]
Note that if any one of these dimensions is finite, then it is equal to those
 on its
left (cf. \cite[Theorem 8.8]{Av}).
\end{remark}

We collect in Theorem \ref{Pneed}
  the  properties of $\Hd_R$
when H is equal to
 Reg (respectively, CI$_*$,  G, or CM) that are used in our proof of  Theorem \ref{lev}.

\begin{theorem}
\label{Pneed}
Let $(R,\m)$ be a local ring and let $M$ be an $R$-module. For H equal
to Reg (respectively, CI$_*$,  G, or CM) the homological dimension $\Hd_R M$ has the
following properties:
\begin{enumerate}[\rm 1.]
\item
If $x$ is $M$-regular and if $\Hd_R M < \infty$ then  $\Hd_R M/xM  < \infty$.

\item
Let $0 \rt N \rt F \rt M \rt 0$ is an exact sequence with $F$ a  free $R$-module.
If $\Hd_R M < \infty$ then  $\Hd_R N < \infty$.

\item
If $N$ is a direct summand of $M$ and if
  $\Hd_R M < \infty$ then  $\Hd_R N < \infty$.

\item
Set $R' = R[X]_{\m R[X]}$ and $M' = M \otimes R'$. If
$\Hd_R M < \infty$ then  $\Hd_R M' < \infty$.
\end{enumerate}
\end{theorem}
\begin{proof}
We first state the results which were known before (See \cite[8.7]{Av}).
 Properties 1  and 2 are  known for each of the homological
 dimensions Reg, CI$_*$,  G, and CM. Property
 3 was known for H = Reg and H = G.  Finally property 4 was known  for
H equal to Reg, G and CI$_*$.
 Proposition \ref{3.3} proves
property 3 in the case H = CI$_*$ and CM. Proposition \ref{3.1} proves  property 4 in
the case H = CM.
\end{proof}

\begin{proposition} \label{3.3} Let $M$ be an $R$-module and let $N$ be a direct summand of $M$.
Then  for H equal
to Reg (respectively, CI$_*$,  G, or CM)  we have $\Hd_RN \leq \Hd_RM.$
\end{proposition}

\begin{proof} This result is known for H= Reg, G. Consider the case
H=CM. We may assume that $\Ct_RM$ is finite for otherwise
there is nothing to prove. Let $R\rightarrow R'\leftarrow Q$ be
the corresponding G-deformation.
 So $\G_Q M\otimes_RR'<
\infty$ and $\Ct_RM=\G_QM\otimes_RR' - \G_QR'$.
Since $N\otimes_RR'$ is a direct summand of $M\otimes_RR'$ we have that
$\G_QN\otimes_RR'\leq \G_Q M\otimes_RR'$.
 In
particular $\Ct_RN$  is finite. Therefore
$\Ct_RN=\G_QN\otimes_RR'-\G_QR'$, So we also get $\Ct_RN  \leq \Ct_RM$.

 Now
suppose H=CI$_*$ and assume that $\CIt M$ is
finite.
Note that for a finite $R$-module $L$ we have
 $$\CIt_R L < \infty  \quad \text{\IFF} \quad \G_R(L) < \infty \ \text{and} \ \cx_R L < \infty. $$
 So $\G_R(M) < \infty \ \text{and} \ \cx_R M < \infty$.
 Since $N$ is a direct summand of $M$ we get that $\G_R(N) < \infty \ \text{and} \ \cx_R N < \infty$. Therefore $\CIt N$ is finite and it is easy to see that
$\CIt N \leq \CIt M $.
\end{proof}

\begin{proposition} \label{3.1} Let $(R,\m)$ be a local ring and $M$ be a finitely
generated $R$-module.
Set $R' = R[X]_{\m R[X]}$ and $M' = M \otimes R'$.
 we have
$\Ct_RM=\Ct_{R'}M'$.
\end{proposition}

To prove this Proposition we need the following Lemma.
\begin{lemma}
\label{flat}
Let $(R,\m)$ and $(S,\n)$ be local rings and let $\phi \colon R \rt S$ be a flat local homomorphism. Set   $R' = R[X]_{\m R[X]}$ and  $S' = S[X]_{\n R[X]}$. Then the naturally
induced $\phi' \colon R' \rt S'$ is a flat local homomorphism.
\end{lemma}
\begin{proof}
Clearly $\phi'$ is a local homomorphism. So it suffices to show $\phi'$ is flat.
Note that $\phi $ induces a ring homomorphism $f \colon R[X] \rt S[X]$. Since
$S$ is a flat $R$-module we have that $S[X] = S\otimes_R R[X]$ is a flat $ R[X]$-module. So $S[X]_{\m R[X]}$ is a flat $R' = R[X]_{\m R[X]}$ module. Since $S'$ is a further localisation of  $S[X]_{\m R[X]}$ we get that $S'$ is a flat $R'$-algebra.
\end{proof}

\begin{proof}[Proof of Lemma \ref{3.1}]
Note that $R'$ is a faithfully flat extension of $R$.
By \cite[8.7(6)]{Av} we have that $\Ct_R M \leq \Ct_{R'} M'$ with equality if
$\Ct_{R'} M'$ is finite.
If $\Ct_RM$
is infinite, then clearly $\Ct_{R'}M'=\infty$ and the
result follows. So suppose $\Ct_RM$ is finite. To show the
equality in this case it is  enough to show that $\Ct_{R'}M'$
is finite. To this end we use definition
3.2 of \cite{Ger}, for CM-dimension. By this definition, if
$\Ct_RM$ is finite, there exists a local flat extension $R
\rightarrow S$ and a suitable $S$-module $K$, such that
$G_K\T{-dim}_S(M \otimes_R S)$ is finite. Since $K$ is suitable, we
have $\T{Hom}_S(K,K)\cong S$ and $\T{Ext}^i_S(K,K)=0$, for all
$i>0.$

Consider the faithfully flat extension $S \rightarrow S'=S[X]_T$
where $T = S[X]\setminus \n S[X]$ and $X$ is an indeterminate over
$S$. Since $S'$ is a flat $S$-module, we
have, $\T{Hom}_{S'}(K \otimes_SS', K \otimes_SS') \cong S'$ and
$\T{Ext}^i_{S'}(K \otimes_SS', K \otimes_SS')=0$ for all $i>0$.
So $K'=K \otimes_SS'$ is a suitable $S'$-module. Moreover using
similar isomorphism, we can deduce that if $\T{G}_K\T{-dim}_SP=0$
then $\T{G}_{K'}\T{-dim}_{S'}P\otimes_SS'=0$. So
$\T{G}_{K'}\T{-dim}_{S'}M\otimes_RS'< \infty$. This in conjunction
with the fact that  $S'$ is a local flat extension of $R'$ (See Lemma \ref{flat}) implies
that $\Ct_{R'}(M \otimes_RR')$ is finite.
\end{proof}

\section{Generalized Levin-Vasconcelos Theorem}
We prove the following:
\begin{theorem}
\label{lev2} Let $(R, \m)$ be a Noetherian local ring with residue
field $k$ and let $M$ be a finite $R$-module of positive depth. If
$\Hd_R(\wt{\m^nM})< \infty$ for some $n$, then $R$ has property H.
\end{theorem}
An easy corollary of this theorem is Theorem \ref{lev}.
\begin{proof}[Proof of Theorem \ref{lev}]
The result follows since for $n \geq \rho(M)$ we have that
$\wt{\m^{n}M} = \m^nM$.
\end{proof}
Theorem \ref{lev2} also has the following corollary:
\begin{corollary}
\label{levc} Let $(R, \m)$ be a Noetherian local ring with residue
field $k$ and let $M$ be a finite $R$-module such that $\depth
G(M)> 0$. If $\Hd(\m^nM)< \infty$ for some $n$, then $R$ has
property H.
\end{corollary}
\begin{proof}
Note that $\depth G(M) > 0 $  implies $\depth M > 0$. Furthermore
in this case $\wt{\m^nM} = \m^n M$ for all $n$. Therefore the
result follows from Theorem \ref{lev2}.
\end{proof}
\begin{proof}[Proof of Theorem \ref{lev2}]
As stated in Remark \ref{HremL}  for the proof of theorem  it is
only enough to consider the cases H= Reg, CI$_*$, G, CM. This we do.

Using the standard trick, in view of Theorem \ref{Pneed}.4  we may assume
that the residue field $R/\m$ is infinite. So there exist $x \in
\m$ which is a $M$-superficial element. Set $N = \wt{\m^n M}$. As
it is shown in Proposition \ref{rrmfull} the submodule $N$ is
$\m$-full and in particular $\m N \colon_M x = N$. Note that $x$
is also $N$-regular. By Theorem \ref{Pneed}.1 we have that  $\Hd_R(N/xN)$ is finite.
 By Proposition
\ref{mfullp}(5) we have that $k$ is a direct summand of $N/xN$. So
by Theorem \ref{Pneed}.3 we get that $\Hd k $ is finite. So $R$ has
property $H$.
\end{proof}
Using Theorem \ref{lev} and Theorem \ref{Pneed}.2 we have the following corollary.
\begin{corollary}
\label{Lcor}
Let $(R,\m)$ be a local ring with $\depth R > 0$.
Let H equal to  Reg (respectively, CI,CI$_*$, G$^*$,G, or CM).
 If
$\Hd_R(R/\m^n)< \infty$ for some $n \geq \rho(R)$, then $R$ has
property H.
\end{corollary}

\section*{Acknowledgement}
We wish to thank L.Avramov and Oana Veliche for carefully reading an earlier draft of this paper and for comments which have improved our exposition. The first
author would like to thank School of Mathematics , TIFR for its hospitality.


\end{document}